\documentclass[11pt,a4paper,twoside,english,leqno]{article}
\usepackage{amsmath,color,amsfonts,amssymb,amscd,amsthm,epsfig}

\usepackage[active]{srcltx}

\usepackage{eqnarray}
\usepackage{enumerate}
\usepackage[english]{babel}
\usepackage[latin1]{inputenc}

\newtheorem{teo}{Theorem}[section]

\newtheorem{prs}[teo]{Proposition}
\newtheorem{prp}[teo]{Property}
\newtheorem{cor}[teo]{Corollary}

\theoremstyle{definition}
\newtheorem{defi}[teo]{Definition}
\theoremstyle{remark}
\newtheorem{ese}[teo]{\textbf{Example}}

\newtheorem{oss}[teo]{ \textbf{Observation}}

\newcommand{\mc}[1]{\mathcal{#1}}
\newcommand{\mb}[1]{\mathbb{#1}}
\newcommand{\be}{\begin{equation}}
\newcommand{\ee}{\end{equation}}
\newcommand{\boss}{\begin{oss}}
\newcommand{\eoss}{\end{oss}}
\newcommand{\bd}{\begin{proof}[\textbf{Proof}]}
\newcommand{\ed}{\end{proof}}
\newcommand{\benonum}{\begin{equation*}}
\newcommand{\eenonum}{\end{equation*}}
\def\è{\grave{e}}
\def\à{\grave{a}}
\def\ò{\grave{o}}
\def\ì{\grave{i}}
\def\ù{\grave{u}}
\def\à{\`a}
\def\è{\`e}
\def\ì{\`\i}
\def\ò{\`o}
\def\ù{\`u}
\def\é{\'e}
\def\({\left\(}
\def\){\right\)}
\def\[{\left\[}
\def\]{\right\]}

\def\ms{\medskip}
\def\noin{\noindent}
\def\sm{\smallskip}

\def\bs{\bigskip}

\newcommand{\ff}{{\mathbb F}}

\newcommand{\zz}{{\mathbb Z}} 

\date{}
\title{Catalan Moments}
\author{\medskip Stefano Barbero, Umberto Cerruti\\ \itshape{Department of Mathematics} \\ \itshape{Università di Torino}\\ \itshape{via Carlo Alberto 8/10 Turin, Italy}}
\begin{document}
\maketitle

\begin{abstract}

This paper is essentially devoted to the study of some interesting relations among the well known operators $I^{(x)}$ (the 
interpolated Invert), $L^{(x)}$ (the interpolated Binomial) and Revert (that we call $\eta$).

We prove that $I^{(x)}$ and $L^{(x)}$ are conjugated in the group $\Upsilon(R)$.\\Here $R$ is a commutative unitary ring. In the same group we see that $\eta$ 
transforms $I^{(x)}$ in $L^{(-x)}$ by conjugation. These facts are proved as corollaries of much more general results.

Then we carefully analyze the action of these operators on the set $\mc{R}$ of second order linear recurrent sequences. 
While $I^{(x)}$ and $L^{(x)}$ transform $\mc{R}$ in itself, $\eta$ sends $\mc{R}$ in the set of moment sequences $\mu_n(h,k)$ of  particular 
families of orthogonal polynomials, whose weight functions are explicitly computed.

The moments come out to be generalized Motzkin numbers (if $R=\zz$, the Motzkin numbers are $\mu_n(-1,1)$). 
We give several interesting expressions of $\mu_n(h,k)$ in closed forms, and one recurrence relation.

There is a fundamental sequence of moments, that generates all the other ones, $\mu_n(0,k)$. These moments are strongly 
related with Catalan numbers. This fact allows us to find, in the final part, a new identity on Catalan 
numbers by using orthogonality relations.
\end{abstract}

\section{A group acting on sequences}

\begin{defi} \label{def1}
\itshape
$$\mathcal{S}(R)=\left\{A=\left\{a_n \right\}_{n=0}^{+\infty}: \ \forall n  \ a_n \in R ,\, a_0 =1\right\}$$

\noin where $R$ is a commutative unitary ring.
\end{defi}

\noin
Now we embed $\mathcal{S}(R)$ in $R[[t]]$ in this way

\be \label{lambd}
\forall A \in \mathcal{S}(R) \quad \quad \lambda(A)=\sum\limits_{n = 0}^{ + \infty } {a_n t^{n + 1} }.
\ee

\noin
In $R[[t]]$ is naturally defined the series composition $\circ$

$$\sum\limits_{n = 0}^{ + \infty } {a_n t^{n + 1} } \circ \sum\limits_{k = 0}^{ + \infty } {b_k t^{k + 1} } = \sum\limits_{n = 0}^{ + \infty } {a_n\left(\sum\limits_{k = 0}^{ + \infty } {b_k t^{k + 1} }\right) ^{n + 1} }.$$

\noin
Then we may induce the operation $\bullet$ in $\mathcal{S}(R)$ :

\begin{defi} \label{def2}

$$\forall A,B\in \mathcal{S}(R) \quad \quad \ A\bullet B=\lambda^{-1}(\lambda(A)\circ \lambda(B))\quad.$$
\end{defi}

\noin
Of course $(\mathcal{S}(R),\, \bullet )$ is a group.
\boss If $R=\ff_q$, this is the \emph{Nottingham group} over $\ff_q$ \cite{Cam} .\eoss

\noin
From every element $A \in \mathcal{S}(R)$ two operators rise: the left multiplication $\mathcal{L}_A$ 
and the right multiplication $\mathcal{R}_A$ .

\begin{defi} \label{def3}

$$\forall B \in \mathcal{S}(R) \quad \quad \mathcal{L}_A(B) = A \bullet B$$

$$\forall B \in \mathcal{S}(R) \quad \quad \mathcal{R}_A(B) = B \bullet A \quad.$$
\end{defi}

\noin We also consider the following two special operators: $\eta$ , often called \textbf{Revert}, and $\varepsilon$ the alternating sign operator :

\begin{defi} \label{def4}

$$ \forall B \in \mathcal{S}(R) \quad \quad \eta(B) = B^{-1} $$

$$\forall B = \left\{b_n \right\}_{n=0}^{+\infty} \in \mathcal{S}(R)  \quad \quad \varepsilon(B) = \left\{(-1)^n b_n \right\}_{n=0}^{+\infty}.$$
\end{defi}

\noin Plainly 

\begin{prp}
\be \label{eta}
\forall A, B \in \mc{S}(R) \quad \quad \eta(A \bullet B) = \eta(B) \bullet \eta(A)
\ee

\be \label{eps}
\forall A, B \in \mc{S}(R) \quad \quad \varepsilon(A \bullet B) = \varepsilon(A) \bullet \varepsilon(B) \quad.
\ee
\end{prp}

\noin In other words, the inversion $\eta$ is an anti-isomorphism of $\mc{S}(R)$, and the alternating sign $\varepsilon$ is an 
isomorphism of $\mc{S}(R)$.

\boss
\noin The operator $\eta$ is especially important.If $a= \left\{a_n\right\}_{n=0}^{+\infty}$ ,  $b = \left\{b_n\right\}_{n=0}^{+\infty}$ and $\eta(a)=b$, then we have 
the relations 

\be \label{e0}
\begin{cases}
u=u(t)=\sum\limits_{n=0}^{+\infty}a_nt^{n+1}\\
t=t(u)=\sum\limits_{n=0}^{+\infty}b_nu^{n+1} \ &\text{inverse series of $u$ .}\\
\end{cases}
\ee
\eoss
\noin 
The operators $\mathcal{L}_A$, $\mathcal{R}_A$, $\eta$, $\varepsilon$, are invertible and can be compounded 
by applying one after the other 
(by the usual operation $\circ$). They generate a group, that we call $\Upsilon(R)$. The group $\Upsilon(R)$ acts on 
$\mathcal{S}(R)$.

\begin{prp} \label{prop1}

For every ring $R$ the following are true

\item \be \label{e1}
\eta = \eta^{-1}
\ee
\item \be \label{e2}
\varepsilon = \varepsilon^{-1}
\ee

\item \be \label{e3}
\eta \circ \varepsilon=\varepsilon \circ \eta
\ee

\item \be \label{e4}
\forall A, B \in \mathcal{S}(R) \quad \quad \mathcal{L}_A \circ \mathcal{R}_B = \mathcal{R}_B \circ \mathcal{L}_A
\ee

\item \be \label{e5}
\forall A \in \mathcal{S}(R) \quad \quad \mathcal{L}_A \circ \eta = \eta \circ \mathcal{R}_{A^{-1}}
\ee

\item \be \label{comee5}
\forall A \in \mathcal{S}(R) \quad \quad  \eta \circ \mathcal{L}_A  =  \mathcal{R}_{A^{-1}} \circ \eta \quad.
\ee

\end{prp}
\bd 
\ms

\noin\\(\ref{e1}) and (\ref{e2}) follow from definition.

\noin\\ (\ref{e3}) Let $d=\eta(\varepsilon(a))$. Because $ \varepsilon(a)=\left\{(-1)^na_n \right\}_{n=0}^{+\infty}$ ,  
(\ref{e0}) becomes
\be 
\begin{cases} \label{e6}
u=u(t)=\sum\limits_{n=0}^{+\infty}(-1)^na_nt^{n+1} \\ 
t=t(u)=\sum\limits_{n=0}^{+\infty}d_nu^{n+1}. \\
\end{cases}
\ee

\noin But (\ref{e6}(\emph{i})) can be rewritten as
$$-u=\sum\limits_{n=0}^{+\infty}a_n(-t)^{n+1}$$
and if $b=\eta(a)$ then
$$t=\sum\limits_{n=0}^{+\infty}(-1)^nb_nu^{n+1}$$ \\
comparing this result with (\ref{e6}(\emph{ii})) we obtain $d=\varepsilon(b)$ .

\noin\\ (\ref{e4})\\
$\quad\quad\quad \forall A,B,C \in \mathcal{S}(R)$
\small
$$(\mathcal{L}_A \circ \mathcal{R}_B)(C) = \mathcal{L}_A(C \bullet B) = A \bullet C \bullet B = \mathcal{R}_B(A \bullet C) = \mathcal{R}_B(\mathcal{L}_A(C)) = (\mathcal{R}_B \circ \mathcal{L}_A)(C)\quad .$$
\normalsize

\noin (\ref{e5}) 
$$ \forall A,B \in \mathcal{S}(R) \quad \quad (\mathcal{L}_A \circ \eta)(B) = A \bullet B^{-1}$$
and
$$ \forall A,B \in \mathcal{S}(R) \quad \quad (\eta \circ \mathcal{R}_{A^{-1}})(B) = (B \bullet A^{-1})^{-1} = A \bullet B^{-1}.$$

\noin (\ref{comee5}) \quad To prove this we do operator composition with $\eta$ both in the front and the back of each side of (\ref{e5})  .
\ed

\bs

\noin Let us pose $\gamma = \eta \circ \varepsilon$ and $X(x) = \left\{x^n \right\}_{n=0}^{+\infty}$, with $x \in R$.

\noin
Of course $X(x) \in \mathcal{S}(R)$ and both $\mathcal{L}_{X(x)}$ and $\mathcal{R}_{X(x)}$ are in 
$\Upsilon(R)$.

\medskip
\noin Plainly

\be \label{lambdaxx} 
\begin{cases}
\text{the generating function of $X(x)$ is } \frac{1}{1-xt}\\
\lambda(X(x)) =  \sum\limits_{n = 0}^{ + \infty } {x^n t^{n + 1} } = \frac{t}{1-x t}.\\
\end{cases}
\ee

\ms

\noin We have

\begin{prp}

\be \label{eps1}
\mc{L}_{X(x)} \circ \varepsilon = \varepsilon \circ \mc{L}_{X(-x)}
\ee

\be \label{eps2}
\mc{R}_{X(x)} \circ \varepsilon = \varepsilon \circ \mc{R}_{X(-x)}
\ee

\be \label{gammainv}
\gamma = \gamma^{-1}
\ee

\be \label{conj1}
\gamma \circ \mc{L}_{X(x)} \circ \gamma^{-1} = \mc{R}_{X(x)}
\ee

\be \label{conj2}
\gamma \circ \mc{R}_{X(x)} \circ \gamma^{-1} = \mc{L}_{X(x)}.
\ee

\end{prp}
\bd 
\bigskip

\noin \\ (\ref{eps1})
 
$\quad \quad \forall A \in \mc{S}(R)$
\small
$$(\mc{L}_{X(x)} \circ \varepsilon)(A) = X(x) \bullet \varepsilon(A) = 
\varepsilon(X(-x)) \bullet \varepsilon(A) =\varepsilon(X(-x) \bullet A) = (\varepsilon \circ \mc{L}_{X(-x)})(A) \quad.$$
\normalsize

\noin\\ (\ref{eps2}) Same proof as for (\ref{eps1}).

\noin\\ (\ref{gammainv}) $\gamma$ is the composition of two commuting involutions.

\noin\\(\ref{conj1})

$\quad\quad\forall A \in \mc{S}(R)$
\small
 $$(\gamma \circ  \mc{L}_{X(x)})(A) = (\eta \circ \varepsilon)(X(x) \bullet A) = 
\eta(X(-x) \bullet \varepsilon(A))= \gamma(A) \bullet X(x) = (\mc{R}_{X(x)} \circ \gamma)(A)\quad.$$
\normalsize 

\noin(\ref{conj2}) Same proof as for (\ref{conj1}).

\ed

\medskip

\noin Let us recall the well known operators Invert and Binomial.

\medskip

 \begin{defi} \label{def5}
 \itshape
The operator $I$ maps the sequence $A = \left\{a_n\right\}_{n=0}^{+\infty}$ in\\ $B = \left\{b_n\right\}_{n=0}^{+\infty}$
where $B$ has generating function:
$$\quad \quad \sum_{n=0}^\infty b_nt^n = \frac{{\sum\limits_{n = 0}^{ + \infty }{a_n t^n } }}{1-t\sum_{n=0}^\infty a_nt^n}\quad.$$
\end{defi}
\medskip
\begin{defi} \label{def6}
\itshape
The operator L maps the sequence $A = \left\{a_n\right\}_{n=0}^{+\infty}$ in \\$B = \left\{b_n\right\}_{n=0}^{+\infty}$ 
where  
$$\quad \quad b_n=\sum_{k=0}^n \binom{n}{k}a_k\quad.$$
\end{defi}

\noin These operators can be \emph{iterated}  \cite{Spi-Ste} and \emph{interpolated} \cite{Bac} becoming  $I^{(x)},\ L^{(y)}$ in this way:
\begin{defi} \label{def7}
\itshape
Given $x \in R$ \ $I^{(x)}$ is called \textbf{Invert} interpolated operator.\ 
By definition  $I^{(x)}(A) = P =\left\{p_n(x)\right\}_{n=0}^{ +\infty }$ where $P$ is the sequence having  generating function

\be \label{invertgen}
\mathbf{P}(t) = \sum\limits_{n = 0}^{ + \infty } {p_n(x) t^n  = \frac{{\sum\limits_{n = 0}^{ + \infty } {a_n t^n } }}{{1 - xt\sum\limits_{n = 0}^{ + \infty } {a_n t^n } }}}\quad.  
\ee 
\end{defi}

\begin{defi}
\itshape Given $y \in R$ \ $L^{(y)}$ is called \textbf{Binomial} interpolated operator. 
By definition
\be \label{bininterp}
 L^{(y)}(A)=\left\{l_n=\sum^{n}_{j=0}{\binom{n}{j} }y^{n - j} a_j \right\}^{+\infty}_{n=0}.
\ee
The exponential generating function of $l=\left\{l_n\right\}_{n=0}^{+\infty}$ is: 
\be \label{binomialgen}
\mathcal{L}(t) = \sum\limits_{n = 0}^{ + \infty }l_n\frac{t^n }{n!} = \sum\limits_{n = 0}^{ + \infty } \sum\limits_{j = 0}^n \frac{(yt)^{n - j} }{(n - j)!}\frac{a_j t^j }{j!} =\exp(ty) \mathcal{A}(t) 
\ee 
being
$$ \exp(ty)=\sum\limits_{n=0}^{+ \infty}\frac{(ty)^n}{n!} \quad \quad \mathcal{A}(t)=\sum\limits_{n=0}^{+\infty}
\frac{a_n t^n}{n!}$$\\
so that (recalling that $a_0=1$) we have  the ordinary generating function
\be \label{binomialser}
\mathbf{L}(t) =\frac{1}{t}A\left( {\frac{t}{{1 - ty}}} \right)
\ee

\noin with $A(t)=\sum\limits_{n=0}^{+\infty}a_nt^{n+1}\quad.$
\end{defi}

\noin The following facts are immediate consequences of (\ref{lambdaxx}).

\begin{prp} \label{xx}
\be \label{xx1}
\forall x \in R \quad \quad \eta(X(x)) = X^{-1}(x) = \left\{(-x)^n \right\}_{n=0}^{+\infty} = X(-x) = \varepsilon(X(x))
\ee
\be \label{xx2}
\forall x, y \in R \quad \quad X(x) \bullet X(y) = X(x+y)\quad.
\ee
\end{prp}
\noin

\noin From their  definitions it is not apparent that the operators $I^{(x)}$ and $L^{(x)}$ are strongly related. Indeed we are 
going to prove that they are, respectively, the left and the right multiplication by $X(x)$ in the group $\mathcal{S}(R)$.

\begin{teo} \label{teo1}

\be \label{teo-i}
 I^{(x)} = \mathcal{L}_{X(x)}
\ee

\be \label{teo-l}
L^{(x)} = \mc{R}_{X(x)}\quad.
\ee

\end{teo}
\bd  
\noin \\(\ref{teo-i}) Let $B=X(x) \bullet A $, then
\footnotesize
$$\lambda(B)=\lambda(X(x))\circ\lambda(A)=
\sum\limits_{n = 0}^{ + \infty } x^n \left( \sum\limits_{k = 0}^{ + \infty } a_k t^{k + 1} \right)^{n+1}= \sum\limits_{k = 0}^{ + \infty } a_k t^{k + 1}  \left(\sum\limits_{n = 0}^{ + \infty } \left( x\sum\limits_{k = 0}^{ + \infty } a_k t^{k + 1} \right)^n \right)=\frac{{\sum\limits_{k = 0}^{ + \infty } {a_k t^{k + 1} } }}{{1 - x\sum\limits_{k = 0}^{ + \infty } {a_k t^{k + 1} } }}$$
\normalsize

\noin so $B=I^{(x)}(A)$ from (\ref{invertgen}) and (\ref{lambd}).

\sm
\noin\\ (\ref{teo-l}) Let $C=A\bullet X(x)$,then
$$\lambda(C)=\lambda(A)\circ\lambda(X(x))=\sum\limits_{n = 0}^{ + \infty } a_n \left(\sum\limits_{k = 0}^{ + \infty } x^k t^{k + 1} \right)^{n+1} = \sum\limits_{n = 0}^{ + \infty } {a_n } \left(\frac{t}{{1 - xt}}\right)^{n + 1}$$ 
so $C=L^{(x)}(A)$ from (\ref{binomialser}) and (\ref{lambd}).

\ed

\noin From Theorem \ref{teo1} and the previous properties we obtain

\begin{teo} \label{teo2}
Let $Id$ be the identity operator and $x,y \in R$.
For the interpolated Invert and Binomial operators the following are true:

$$I^{(x)} \circ I^{(-x)} = Id \quad\quad L^{(x)} \circ L^{(-x)} = Id$$

$$I^{(x)} \circ I^{(y)} = I^{(x+y)} \quad\quad L^{(x)} \circ L^{(y)} = L^{(x+y)}$$

$$I^{(x)} \circ \varepsilon = \varepsilon \circ I^{(-x)} \quad\quad L^{(x)} \circ \varepsilon = \varepsilon \circ L^{(-x)}$$

$$I^{(x)} \circ \ L^{(y)} = L^{(y)} \circ I^{(x)}$$

$$I^{(x)} \circ \eta = \eta \circ L^{(-x)} \quad \quad  \eta \circ I^{(x)} = L^{(-x)} \circ \eta$$

$$\gamma \circ I^{(x)} \circ \gamma^{-1} = L^{(x)} \quad \quad \gamma \circ L^{(x)} \circ \gamma^{-1} = I^{(x)}\quad.$$

\end{teo}

\medskip

\noin So we have seen, by the way, that the operators $I^{(x)}$ and $L^{(x)}$ are \emph{conjugated} in the group $\Upsilon(R)$!

\section{Action on linear recurrent sequences of order 2}

\noin 
In this section we analyze the action of $I^{(x)}$ and $L^{(x)}$ on the particular subset of $\mc{S}(R)$ 
formed by linear recurrent sequences of order 2 (starting with 1).
\begin{defi}
\itshape
$$\mathcal{R}(R)=\left\{\mathcal{W}(1,b,h,k):b,h,k \in R\right\}$$\\
where $$\mathcal{W}(1,b,h,k)=\left\{\mathcal{W}_n(1,b,h,k)\right\}_{n=0}^{+\infty}$$\\
satisfies the recurrence $\forall n \geq 2$
\be \label{rec}
\left\{ \begin{array}{l}
 {\mathcal{W}_0(1,b,h,k)  = 1} \\ 
 {\mathcal{W}_1(1,b,h,k)  = b} \\ 
 \mathcal{W}_n(1,b,h,k)  = h\mathcal{W}_{n-1}(1,b,h,k)  - k\mathcal{W}_{n-2}(1,b,h,k)\quad \forall n \geq 2\quad.\\ 
 \end{array} \right.
\ee
\end{defi}

\noin $I^{(x)}$ and $L^{(x)}$ map $\mathcal{R}(R)$ into itself in the following way

\begin{teo}\label{recteo}
\itshape

$\forall x,y \in R$ we have
\be \label{A} I^{(x)}(\mathcal{W}(1,b,h,k))=\mathcal{W}(1,b+x,h+x,(h-b)x+k)\ee
\be \label{B} L^{(y)}(\mathcal{W}(1,b,h,k))=\mathcal{W}(1,b+y,h+2y,y^2+hy+k)\ee
\small
\be \label{C}\\ C^{(x,y)}(\mathcal{W}(1,b,h,k))=\mathcal{W}(1,b+y+x,h+x+2y,y^2+hy+k+(h-b)x+xy)\ee
\normalsize
where $C^{(x,y)}=I^{(x)}\circ L^{(y)}=L^{(y)} \circ I^{(x)}$.

\end{teo}
\bd

The generating function of $\mathcal{W}(1,b,h,k))$ is 
\be \label{genf}\mathbf{W}(t) = \frac{{1 + (b - h)t }}{{1 - ht + kt^2 }}.\ee\\
\ms
If we substitute $\mathbf{W}(t)$ in (\ref{invertgen}), and compute $\mathbf{P}(\mathbf{W}(t))$, we find 
$$\frac{1+(b-h)t}{1 -(h+x)t +(k +(h-b)x)t^2}.$$ This proves (\ref{A}).
\medskip
\noin In the ring $R[z]/(z^2-hz+k)$ we pose $\alpha_1=z$ and $\alpha_2=h-z$ (the roots of $z^2-hz+k$). Then we have 
$\mathcal{W}_n(1,b,h,k)) = p\alpha_1^n+q\alpha_2^n$.

\noin \\ Substituting the sequence $\mathcal{W}(1,b,h,k))$ to the sequence $A$ in (\ref{bininterp}) we obtain 
\small
$$ l_n  = \sum\limits_{i = 0}^{n} {\binom{n}{i}} y^{n -i} a_{i}  = \sum\limits_{i = 0}^{n} {\binom{n}{i}}y^{n - i} (p\alpha _1^{i}  + q\alpha _2^{i} ) =$$$$=p\sum\limits_{i = 0}^{n} {\binom{n}{i}} y^{n - i} \alpha _1^i  + q\sum\limits_{i = 0}^{n} {\binom{n}{i}} y^{n  - i} \alpha _2^i=p(y + \alpha _1 )^{n}  + q(y + \alpha _2 )^{n}\quad.$$\\
\normalsize
Then posing $y + \alpha _1=R$ and $y + \alpha _2=S$, observing that
\ms
\small
$$R+S=2y+\alpha _1  + \alpha _2 =h + 2y \quad RS=y^2  + (\alpha _1  + \alpha _2 )y + \alpha _1 \alpha _2 =y^2  + hy + k$$
\normalsize
\ms
where we used $\alpha _1  + \alpha _2=h$ and $\alpha _1 \alpha _2=k$, we have:
\small
$$ 
 l_n= pR^{n}  + qS^{n}  + pR^{n-1}S + qS^{n-1}R -pR^{n-1}S - qS^{n-1}R = $$  
 $$ = R(pR^{n - 1} + qS^{n - 1})+S(pR^{n - 1} + qS^{n - 1})-RS(pR^{n - 2}  + qS^{n - 2} ) = $$ 
 $$ = (R+S)(pR^{n - 1}  + qS^{n - 1} )- RS(pR^{n - 2}  + qS^{n - 2} )= (h + 2y)l_{n - 1}  - (y^2  + hy + k)l_{n - 2}\quad.$$ 
\normalsize
 \noin This proves (\ref{B}).
 \medskip
 
 \noin (\ref{C}) follows at once from (\ref{A}) and (\ref{B}).
\ed

\noin \boss
An important subset $\mathcal{F}\subset\mathcal{R}(R)$ consists of sequences $$F(h,k)=\mathcal{W}(1,h,h,k)=\left\{1,h,h^2-k,..\right\}$$
\be 
\left\{ \begin{array}{l}
{\mathcal{W}_0(1,h,h,k)  = 1} \\ 
{\mathcal{W}_1(1,h,h,k) = h }\\ 
{\mathcal{W}_n(1,h,h,k)  = h\mathcal{W}_{n-1}(1,h,h,k)  - k\mathcal{W}_{n-2}(1,h,h,k) \quad \forall n \ge 2}\quad.\\ 
 \end{array} \right. \ee
They are a subset of generalized Fibonacci sequences.\\From Theorem \ref{recteo} we can define a polynomial sequence $\mathcal{P}(h,k,x)$ as follows
\be \label{recfib}\mathcal{P}(h,k,x)=\left\{P_n(h,k,x)\right\}_{n=0}^{+\infty}=I^{(x)}(F(h,k))=F(h+x,k)\ee
and we can observe that
 \be I^{(h)}(\mathcal{W}(1,0,0,k))=I^{(h)}(F(0,k))=F(h,k).\ee
These relations, as we will see in the next section, show a connection between $\mathcal{F}$ and orthogonal polynomials. They also help us to prove  what Bacher \cite{Bac} observes about arithmetical properties of $P_n(h,k,x)$.
 \eoss
 \begin{prs}
 \itshape
  $ \forall$ $ m,n$ such that $m|n$ then $P_{m-1}(h,k,x)|P_{n-1}(h,k,x)$. 
 \end{prs}
 \bd
 \noin \\(\ref{recfib}) gives the recurrence relation
 \be \label{recpol}
\left\{ \begin{array}{l}
 P_0(h,k,x)  = 1 \\ 
 P_1(h,k,x)  = h+x \\ 
 P_n(h,k,x)  = (h+x)P_{n - 1}(h,k,x)  - kP_{n - 2}(h,k,x) {\rm   } \ \forall n \ge 2 \\ 
 \end{array} \right.
\ee
 from which 
 $$ P_n(h,k,x)  = \frac{{\alpha _1^{n+1}  - \alpha _2^{n+1} }}{{\alpha _1  - \alpha _2 }} $$\\
 where $\alpha _1$ and $\alpha _2$ are the roots of the characteristic polynomial 
 $$ t^2-(h+x)t+k=0 \quad.$$ \\
Thus
$$ \frac {{P_{n-1}(h,k,x)}}{{P_{m-1}(h,k,x)}}  = \frac{{\alpha _1^{n}  - \alpha _2^{n} }}{{\alpha _1^{m}  - \alpha _2^{m} }}$$\\
and if $m|n$ then $P_{m-1}(h,k,x)|P_{n-1}(h,k,x)$. 
 
 \ed
 \noin Finally we can find a  couple  of  relations  on  sequences  $\mathcal{W}(1,b,h,k)$ \\
 involving the $\eta$ operator.
\begin{cor}\label{core}
\itshape
For all sequences $\mathcal{W}(1,b,h,k)$ we have

\be \label{S} I^{(x)}(\eta(\mathcal{W}(1,b,h,k)))=\eta(W(1,b-x,h-2x,x^2-hx+k))\ee
\be \label{T} L^{(x)}(\eta(\mathcal{W}(1,b,h,k)))=\eta(W(1,b-x,h-x,(b-h)x+k))\quad.\ee
\end{cor}
\bd
The proof is obvious from Theorem (\ref{teo2}).
\ed

\section{Moments generating function}
From now on we shall pose $R=\mb{C}$.

\medskip

\noin We know from (\ref{recfib}) that $I^{(x)}$, applied to elements in $\mathcal{F}$, gives rise to polynomial sequence $\mathcal{P}(h,k,x)=\left\{P_n(h,k,x)\right\}_{n=-1}^{+\infty}$ (where indexes have been changed in (\ref{recfib}) for convenience in calculation), with recurrence relation
\be \label{pol}
\left\{ \begin{array}{l}
 P_{-1} (h,k,x) = 0 \\ 
 P_0 (h,k,x)  = 1 \\ 
 P_n (h,k,x) = (x + h)P_{n - 1} (h,k,x) - kP_{n - 2} (h,k,x)\quad  \forall n\geq1 \quad. \\ 
 \end{array} \right.
 \ee
 \\ From  Favard's theorem (\cite{Chi}, page 21)
  this recurrence relation, when $k\neq0$, is also the one for orthogonal polynomials having a proper moment functional. If $h=0$ we have $P_n(0,k,x)=E_n(x,k)$ the $n$-th Dickson polynomial of the second kind \cite{Lidl}.\\ Moreover for the moments sequence $\mu(h,k)$ related to the sequence $\mathcal{P}(h,k,x)$ the following  holds : 
\begin{teo}
\itshape
The sequence $\mu(h,k)$ has generating function
\be \label{mgf}
\mu(t)=\sum\limits_{n = 0}^{ + \infty }\mu _n t^n  = \frac{{1 - ht - \sqrt {(1 - ht)^2  - 4kt^2 } }}{{2kt^2 }} \quad.
\ee
\end{teo}
\bd
From known results about orthogonal polynomials theory \cite{Chi},\ the moments generating function $\mu(t)$ is equal to a continued fraction :
\be\label{contfr}\mu(t)=\sum\limits_{n = 0}^{ + \infty } {\mu _n t^n }  = \frac{{\lambda _0 }}{{1 + \xi _0 t - \frac{{\lambda _1 t^2 }}{{1 + \xi _1 t - \frac{{\lambda _2 t^2 }}{{1 + \xi _2 t - \frac{{\lambda _3 t^2 }}{{1 + \xi _3 t - .....}}}}}}}} \quad.\ee
For $\mathcal{P}(h,k,x)$ we have $\forall n$ $\xi_n=-h$,\ $\forall n>1$ $\lambda_n=k$, $\lambda_0=\mu_0=1$ and
 (\ref{contfr}) becomes
\be 
\mu(t)=\sum\limits_{n = 0}^{ + \infty } {\mu _n t^n }  = \frac{1}{{1 - ht - \frac{{kt^2 }}{{1 - ht - \frac{{kt^2 }}{{1 -ht - \frac{{kt^2 }}{{1 - ht - .....}}}}}}}}\quad.\ee
It can be expressed in closed form posing 
$$ 
y = \frac{{kt^2 }}{{1 - ht - \frac{{kt^2 }}{{1 - ht - \frac{{kt^2 }}{{1 - ht - \frac{{kt^2 }}{{1 - ht - .....}}}}}}}}
$$\\ and observing that \be\label{oss} 
y = \frac{{kt^2 }}{{1 - ht - y}}
\ee
and
\be \label{mou}
\mu (t) = \frac{1}{{1 - ht - y}}\quad.\ee \noin Finding $y$ from (\ref{oss}) we obtain 

$$
y_1 = \frac{{1 - ht + \sqrt {(1 - ht)^2  - 4kt^2 } }}{2}$$  $$ y_2 = \frac{{1 - ht - \sqrt {(1 - ht)^2  - 4kt^2 } }}{2}\quad.
$$
We have to choose $y=y_2$ because $y_1$ replaced in (\ref{mou}) gives rise to discontinuity at $t=0$.
With this value for $y$  and a rationalization we easily find the exact form of $\mu(t)$ in  (\ref{mgf}).
\ed
\noin The explicit moments values are given by the
\begin{cor}
\itshape
The moments $\mu_n(h,k)$ related to polynomials $\mathcal{P}(h,k,x)$ are equal to
\be \label{mom}
\mu_n(h,k)=\begin{cases}
 - \frac{1}{{2k}}\sum\limits_{j = 0}^{\frac{{n + 1}}{2}} {\binom{1/2}{n +2-j}}{\binom{n + 2 - j}{j} }( - 2h)^{n + 2 - 2j} (h^2  - 4k)^j  &\text{$n$  odd}\\
- \frac{1}{{2k}}\sum\limits_{j = 0}^{\frac{{n + 2}}{2}} {\binom{1/2}{n +2-j}}{\binom{n + 2 - j}{j} }( - 2h)^{n + 2 - 2j} (h^2  - 4k)^j  &\text{ $n$ even}\\
\end{cases}
\ee
where $n\geq 1$ and $\mu_0=1$.
\end{cor}
\bd
The result follows developing (\ref{mgf}):
$$
\sqrt {(1 - ht)^2  - 4kt^2 }  = (1 - 2ht + (h^2  - 4k)t^2 )^{1/2}  = \sum\limits_{i = 0}^{ + \infty } {\binom{1/2}{i}} (-2ht + (h^2  - 4k)t^2 )^i $$\\
and so
\small
$$
\mu (t) = \frac{1}{{2kt^2 }}\left(1 - ht - (1 - ht) - \frac{1}{2}(h^2  - 4k)t^2  - \sum\limits_{i = 2}^{ + \infty } {\binom{1/2}{i}}(-2ht + (h^2  - 4k)t^2 )^i \right) = $$
$$
 =  - \frac{1}{{4k}}(h^2  - 4k) - \frac{1}{{2k}}\sum\limits_{i = 2}^{ + \infty } {\binom{1/2}{i}} \sum\limits_{j = 0}^i {\binom{i}{j}} (-2h)^{i - j} (h^2  - 4k)^j t^{i + j - 2}\quad.$$\\
\normalsize
Ordering the summation with respect to the degree $n$ of  $t^{n}$ , we observe that the coefficient of $t^{n}$ for $n=0$ is $-h^2/2$ and replacing it in $\mu(t)$ expression we have
$$
\mu(t)=1+\sum\limits_{n=1}^{+\infty}\mu_n(h,k)t^n$$

$$\mu_n(h,k)=\begin{cases}
 - \frac{1}{{2k}}\sum\limits_{j = 0}^{\frac{{n + 1}}{2}} {\binom{1/2}{n +2-j}}{\binom{n + 2 - j}{j}}( - 2h)^{n + 2 - 2j} (h^2  - 4k)^j  &\text{for odd $n$ }\\
- \frac{1}{{2k}}\sum\limits_{j = 0}^{\frac{{n + 2}}{2}} {\binom{1/2}{n +2-j}}{\binom{n + 2 - j}{j} }( - 2h)^{n + 2 - 2j} (h^2  - 4k)^j  &\text{for even $n$ .}\\
\end{cases}$$
\ed
\boss
\noin The moments $\mu_n(h,k)$  are the generalized Motzkin numbers.We will show a combinatorial interpretation of them in Section 6.
\eoss
\section{Weight function}
We want to find the weight function  $\omega(t)$ of the functional $\mathcal{V}$ related to the sequence $\mathcal{P}(h,k,x)$ (see \cite{Chi}). So $\mathcal{V}[f]$ will be defined as follows
$$\mathcal{V}[f]=\int_{C}f(t)d\psi(t)$$\\
where \textsl{C} will be a suitable integration interval, $\psi(t)$ a distribution such that $\psi^{'}(t)=\omega(t)$. By Stieltjes inversion formula we have
\be\label{invfor}
\psi (t) - \psi (0) =  - \frac{1}{\pi }\mathop {\lim }\limits_{y \to 0^ +  } \int\limits_0^t {{\mathop{\mathcal{I}m}\nolimits} (F(x + iy,h,k))dx} 
\ee
being $z=x+iy \in\mb{C}$ and  $F(z,h,k)=z^{-1}\mu(z^{-1})$ where $\mu(t)$ is defined by (\ref{mgf}); thus
\be F(z,h,k)=\frac{{z - h - \sqrt {(z - h)^2  - 4k} }}{{2k}}\quad.\ee
\noin We can immediately find the corresponding primitive $\mathcal{F}(z,h,k)$ of $F(z,h,k)$
\scriptsize
$$\mathcal{F}(z,h,k)= \frac{1}{{2k}} \left(\frac{{z^2 }}{2} - hz - \frac{{(z - h)}}{2}\sqrt {(z - h)^2  - 4k}-2k\log \left(\sqrt {(z - h)^2  - 4k}  - (z - h)\right)\right)
$$\\
\normalsize
where the arbitrary constant has been made equal to 0, without loss of generality.\\
Now we can study, depending on $h$,\ $k$, the value of
\be\label{lim}
{\mathop{\mathcal{I}m}\nolimits} \left(\mathop {\lim }\limits_{y \to 0^ +  } \mathcal{F}(x + iy,h,k)\right)
\ee
considering all the parts which summed together give $\mathcal{F}$ :
\begin{enumerate}
\item \small $$\mathop{\lim}\limits_{y\to 0^+}\left(\frac{1}{2k}\left(\frac{{z^2}}{2}-hz\right)\right)=\frac{1}{2k}\left(\frac{{x^2}}{2}-hx\right)$$ \normalsize
\item \small
$$
\mathop {\lim }\limits_{y \to 0^ +  } \left(\frac{1}{{4k}}(z - h)\sqrt {(z - h)^2  - 4k} \right) = \frac{1}{{4k}}(x - h)\sqrt {(x - h)^2  - 4k} 
$$
\normalsize
\item \small
$$ 
\mathop {\lim }\limits_{y \to 0^ +  } \left(\log \left(\sqrt {(z - h)^2  - 4k}  - (z - h)\right)\right) = \log \left(\sqrt {(x -h)^2  - 4k}  - (x - h)\right)
$$
\normalsize
\end{enumerate}
remembering the condition $k\neq 0$,\ we note that :

\noin \\ i)\ is always real;

\noin \\ ii)\ is real if $(x-h)^2-4k\geq 0$ or $k< 0$,\ otherwise
$$
\mathop {\lim }\limits_{y \to 0^ +  } \left(\frac{1}{{4k}}(z - h)\sqrt {(z - h)^2  - 4k} \right) = \frac{i}{{4k}}(x - h)\sqrt {4k-(x - h)^2} 
$$
when $h-2\sqrt{k}<x<h+2\sqrt{k}$ ;

\noin \\ iii)\ if $k<0$ or $k>0$ and $x\notin( h-2\sqrt{k},h+2\sqrt{k})$,\ $\sqrt{(x-h)^2-4k}$  is real, moreover
$$ \sqrt{(x-h)^2-4k}-(x-h)>0$$ surely if $k<0$,while if $k>0$ and $x\notin( h-2\sqrt{k},h+2\sqrt{k})$ the logarithm is real
 if $x\in(-\infty,h-2\sqrt{k})$ and complex if $x\in( h+2\sqrt{k},+\infty).$\\
In this ultimate case we have 
$$
\log\left(\sqrt{(x-h)^2-4k}-(x-h)\right)=\log\left|\sqrt{(x-h)^2-4k}-(x-h)\right|+i\pi \quad .$$\\
Finally if $k>0$ and $x\in( h-2\sqrt{k},h+2\sqrt{k})$ then $$\sqrt{(x-h)^2-4k}=i\sqrt{4k-(x-h)^2}$$ and
\small 
$$
\log\left(\sqrt{(x-h)^2-4k}-(x-h)\right)=\log\left(-(x-h)+i\sqrt{4k-(x-h)^2}\right)=$$
$$=\log\left|-(x-h)+i\sqrt{4k-(x-h)^2}\right|+iArg\left(-(x-h)+i\sqrt{4k-(x-h)^2}\right)$$\\
\normalsize
with
\small
$$
Arg\left(-(x-h)+i\sqrt{4k-(x-h)^2}\right)=\begin{cases}
 -\arctan \left( {\frac{{\sqrt {4k - (x - h)^2 } }}{{ (x - h)}}} \right)&\text{ if $ h - 2\sqrt k  < x <   h $}\\ 
 \frac{\pi }{2} &\text{if $x =   h$} \\ 
 \pi  - \arctan \left( {\frac{{\sqrt {4k - (x - h)^2 } }}{{(x - h)}}} \right)&\text{ if $ h < x <  h + 2\sqrt k$ .} \\ 
  
 \end{cases}$$\\
 \normalsize
 So the  limit (\ref{lim}) is zero  for $k<0$ and also for $k>0$ with $ x\in(-\infty,h-2\sqrt{k})$ while
 when \ $k>0$ and $x\in( h-2\sqrt{k},+\infty)$ the limit values are
 \small
$$
\begin{cases}
 
  - \frac{{(x - h)\sqrt {4k - (x - h)^2 } }}{{4k}} +\arctan \left( {\frac{{\sqrt {4k - (x - h)^2 } }}{{ (x - h)}}} \right) &\text{           if $ h - 2\sqrt k  < x <   h $}\\ 
  - \frac{\pi }{2}  &\text{ if $x =   h$} \\ 
  - \frac{{(x - h)\sqrt {4k - (x - h)^2 } }}{{4k}} - \pi  + \arctan \left( {\frac{{\sqrt {4k - (x - h)^2 } }}{{(x - h)}}} \right)     &\text{ if $ h < x <   h + 2\sqrt k$}  \\ 
  - \pi &\text{ if $x >   h + 2\sqrt k$  .}  \\ 
 \end{cases} $$\\
 \normalsize 
 This gives, together with (\ref{invfor}) 
 \be\label{omega}
 \omega(t)=\psi'(t)= \begin{cases}
 
 \frac{{\sqrt {4k - (t - h)^2 } }}{{2k\pi }}\   &\text{if  $ h - 2\sqrt k  < t <   h + 2\sqrt k \wedge t \ne   h $}\\ 
 0\ &\text{otherwise .} \\ 
 \end{cases}
  \\ 
\ee
\section{Recurrence relation for  $\mu(h,k)$}
We know from definition that
\be \label{moment}\mu_n=\mathcal{V}[t^n]=\int_{C}t^n d\psi(t)\ee
and relation  (\ref{moment}) becomes, using (\ref{omega})
\be \label{int}\mu_n=\int^{  h + 2\sqrt k }_{  h - 2\sqrt k } {\frac{{t^n \sqrt {4k - (t - h)^2 } }}{{2k\pi }}} dt\quad .\ee
Now we can prove the 
\begin{teo}
\itshape
The sequence $\mu(h,k)$ is  recurrent with
\be
\begin{cases}
\mu_0=1\\
\mu_1=h\\
\mu_n=\frac{{h(2n+1)\mu_{n-1}-(h^2-4k)(n-1)\mu_{n-2}}}{n+2} &\text{$\forall n\geq 2$ .}\\
\end{cases} \ee

\end{teo}
\bd
$ \forall n\geq 2 $ we have
$$
\mu_n\quad=\int^{  h + 2\sqrt k }_{ h - 2\sqrt k } \frac{{t^{n-1}(t-h+h) \sqrt {4k - (t - h)^2 } }}{{2k\pi }} dt=$$
$$ \quad\quad\quad=\int^{ h + 2\sqrt k }_{  h - 2\sqrt k } \frac{{t^{n-1}(t-h) \sqrt {4k - (t - h)^2 } }}{{2k\pi }} dt +h\mu_{n-1}$$\\
using integration by parts we obtain 
 \small
$$\mu_n=
\left[ { \frac{{-t^{n - 1} \sqrt {(4k - (t - h)^2 )^3 } }}{{6k\pi }}} \right]_{h - 2\sqrt k }^{h + 2\sqrt k }  - \int^{h + 2\sqrt k }_{h - 2\sqrt k } {\frac{{ - (n - 1)t^{n - 2} \sqrt {(4k - (t - h)^2 )^3 } }}{{6k\pi }}}dt+h\mu_{n-1}
$$
\normalsize
but 
$$\left[ { \frac{{-t^{n - 1} \sqrt {(4k - (t - h)^2 )^3 } }}{{6k\pi }}} \right]_{h - 2\sqrt k }^{h + 2\sqrt k }=0$$
so 
\small
$$\mu_n=-\frac{{(n-1)}}{3}\int^{  h + 2\sqrt k }_{  h + 2\sqrt k } {\frac{{t^{n-2}(h^2-4k+ t^2-2ht) \sqrt {4k - (t - h)^2 } }}{{2k\pi }}} dt+h\mu_{n-1}=$$
$$\quad\quad =-\frac{{(h^2-4k)(n-1)}}{3}\mu_{n-2}-\frac{{(n-1)}}{3}\mu_n+\frac{{2h(n-1)}}{3}\mu_{n-1}+h\mu_{n-1}=$$
$$=-\frac{{(h^2-4k)(n-1)}}{3}\mu_{n-2}-\frac{{(n-1)}}{3}\mu_n+\frac{{h(2n+1)}}{3}\mu_{n-1}$$\\
\normalsize
we finally find the recurrence 
$$\mu_n=\frac{{h(2n+1)\mu_{n-1}-(h^2-4k)(n-1)\mu_{n-2}}}{n+2}$$\\while  $\mu_0$ and $\mu_1$ can be easily found calculating (\ref{int}) for $n=0,1$ .
\ed
\begin{cor}\label{corL}
\itshape
We have $\forall y \in \mb{R} \ L^{(y)}(\mu(h,k))=\mu(h+y,k)$ .
\end{cor}
\bd
In fact if
$$\mu_n^{'}=\sum\limits_{i=0}^{n}{\binom{n}{i}}y^{n-i}\mu_i$$\\
 using (\ref{int}) we have
$$\mu_n^{'}=\int^{  h + 2\sqrt k }_{  h + 2\sqrt k } {\sum\limits_{i=0}^{n}{\binom{n}{i}}y^{n-i}t^i\frac{ {\sqrt {4k - (t - h)^2 } }}{{2k\pi }}}dt\\$$
Now
$$\sum\limits_{i=0}^{n}{\binom{n}{i}}y^{n-i}t^i=(t+y)^i$$\\
and substituting $u=t+y$ and $h^{'}=h+y$ 
$$\mu_n'=\int^{  h^{'} + 2\sqrt k }_{ h^{'} - 2\sqrt k } {\frac{{u^n \sqrt {4k - (u -h^{'})^2 } }}{{2k\pi }}}du$$\\
Thus $\mu_n^{'}$ is defined with an analogous relation like (\ref{int}) for $\mu_n$ .
\ed
\section{Combinatorial interpretation for $\mu_n(h,k)$ }
We consider a lattice $(n+1)\times (n+1)$ composed by all the points having non negative integer coordinates. Motzkin paths are all
the courses starting from $(0,0)$ and reaching $(n,0)$ with the following rules
$$\begin{cases}
(i,j)\rightarrow(i+1,j) &\text{horizontal shift  to east}\\
(i,j)\rightarrow(i+1,j+1) &\text{diagonal shift to  north-east}\\
(i,j)\rightarrow(i+1,j-1) &\text{diagonal shift to south-east}\\
\end{cases}$$
For example from $(0,0)$  to $(3,0)$  we have only the 4 possible paths
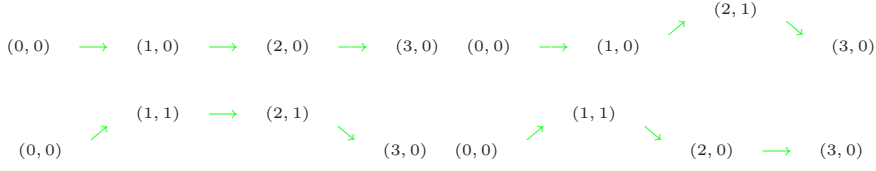
\tiny
$$
\begin{array}{*{20}c}
   {} & {} & {} & {} & {} & {} & {}  \\
   {} & {} & {} & {} & {} & {} & {}  \\
   {(0,0)}  & {\textcolor{green}{\longrightarrow}} &{(1,0)}  & {\textcolor{green}{\longrightarrow}} &{(2,0)} & {\textcolor{green}{\longrightarrow}} &{(3,0)}  \\
\end{array}
\begin{array}{*{20}c}
   {} & {} & {} & {} &  {(2,1)}  & {} & {}  \\
   {} & {} & {} & {\textcolor{green}{\nearrow}} & {} & {\textcolor{green}{\searrow}} & {}  \\
   {(0,0)}  & {\textcolor{green}{\longrightarrow}} &{(1,0)}  & {} & {} & {} &{(3,0)}   \\
\end{array}
$$
\medskip
$$
 \begin{array}{*{20}c}
   {} & {} & {(1,1)} &{\textcolor{green}{\longrightarrow}}  & {(2,1)} & {} & {}  \\
   {} & {\textcolor{green}{ \nearrow}}  & {} & {} & {} &  {\textcolor{green}{\searrow}}  & {}  \\
   {(0,0)} & {} & {} & {} & {} & {} & {(3,0)}  \\
\end{array}\\ \\
\begin{array}{*{20}c}
   {} & {} & {(1,1)} & {} & {} & {} & {}  \\
   {} &  {\textcolor{green}{\nearrow}}  & {} & {\textcolor{green}{ \searrow}}  & {} & {} & {}  \\
   {(0,0)} & {} & {} & {} & {(2,0)} &{\textcolor{green}{\longrightarrow}} & {(3,0)}  \\
\end{array} 
$$
 \normalsize
\noin \\ If we weight one shift of a path $\mathcal{P}$ posing:
 $$\begin{cases}
w((i,j)\rightarrow(i+1,j))=h \\
w((i,j)\rightarrow(i+1,j+1))=1 \\
w((i,j)\rightarrow(i+1,j-1))=k \\
\end{cases}$$
we can describe $\mathcal{P}$ with weights product.
The four paths represented above are respectively represented by:\ $h^3$,\ $hk$,\ $hk$,\ $kh$.\ 
We observe that the sum of all the weights of these paths from (0,0) to (3,0) is $h^3+3hk=\mu_3(h,k)$.\ This is a consequence of the 
\begin{teo}[Viennot's Theorem \cite{Sta}]
\itshape
Under the rules described above,for every Motzkin path $P$ the following relation holds 
$$\mu _n (h,k) = \sum\limits_{P:(0,0) \to (n,0)} {w(P)}\quad. $$
\end{teo}
\noin As a consequence in $\mu_n(h,k)$ is codified  information about all weighted paths from $(0,0)$ to $(n,0)$:
\begin{itemize}
\item the sum of coefficients of $\mu_n(h,k)$ gives the number of all possible Motzkin paths from $(0,0)$ to $(n,0)$;
\item the  $h$ exponent in every term gives the number of horizontal shifts to east;
\item the $k$ exponent in every term gives the number of diagonal shifts to  north-east (or to south-east);
\item the weight $h$ may be interpreted as the number of colors among which we can select one to draw horizontal shifts to east;
\item the weight $k$ may be interpreted as the number of colors among which we can select one to draw diagonal shifts to north-east or to south-east;
\end{itemize}
\begin{ese} \noin \\From $\mu_4(h,k)=h^4+6h^2k+2k^2$ we have $9=1+6+2=\mu_4(1,1)$ distinct paths from $(0,0)$ to $(4,0)$ traced with one color for all shifts:
\begin{itemize}
\item 1 path having 4 horizontal shifts;
\item 6 with 2 horizontal shifts and 1 to north-east(and so 1 to south-east);
\item 2 with 2 diagonal shifts to north-east (and so 2 to south-east).
\end{itemize}
Moreover from $\mu_3(1,1)=4$ we recover the previous one-colored paths and from $\mu_3(1,2)=7$ we find all the paths painted with one color for horizontal shifts  and two possible colors  for diagonal shifts:
\end{ese}
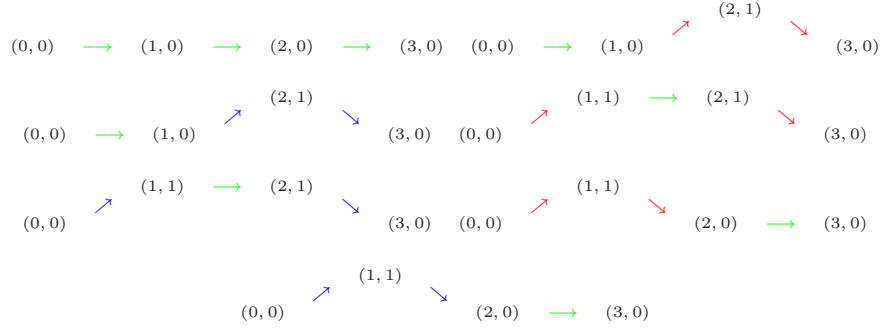
\tiny
$$
\begin{array}{*{20}c}
   {} & {} & {} & {} & {} & {} & {}  \\
   {} & {} & {} & {} & {} & {} & {}  \\
   {(0,0)} & {\textcolor{green}{\longrightarrow}} & {(1,0)} & {\textcolor{green}{\longrightarrow}} & {(2,0)}& {\textcolor{green}{\longrightarrow}} & {(3,0)}   \\
\end{array}\\
\begin{array}{*{20}c}
   {} & {} & {} & {} &{(2,1)} & {} & {}  \\
   {} & {} & {} & {\textcolor{red}{\nearrow}} & {} & {\textcolor{red}{\searrow}} & {}  \\
   {(0,0)} & {\textcolor{green}{\longrightarrow}} & {(1,0)} & {} & {} & {} & {(3,0)}   \\
\end{array}
$$
\medskip
$$
\begin{array}{*{20}c}
   {} & {} & {} & {} & {(2,1)} & {} & {}  \\
   {} & {} & {} & {\textcolor{blue}{\nearrow}} & {} & {\textcolor{blue}{\searrow}} & {}  \\
   {(0,0)} & {\textcolor{green}{\longrightarrow}} & {(1,0)} & {} & {} & {} & {(3,0)}   \\
\end{array}
\begin{array}{*{20}c}
   {} & {} & {(1,1)} & {\textcolor{green}{\longrightarrow}} & {(2,1)} & {} & {}  \\
   {} & {\textcolor{red}{\nearrow}} & {} & {} & {} & {\textcolor{red}{\searrow}} & {}  \\
   {(0,0)} & {} & {} & {} & {} & {} & {(3,0)}  \\
\end{array}
$$
\medskip
$$
\begin{array}{*{20}c}
   {} & {} &  {(1,1)} & {\textcolor{green}{\longrightarrow}} & {(2,1)}  & {} & {}  \\
   {} & {\textcolor{blue}{\nearrow}} & {} & {} & {} & {\textcolor{blue}{\searrow}} & {}  \\
   {(0,0)}  & {} & {} & {} & {} & {} &  {(3,0)}  \\
\end{array}
\begin{array}{*{20}c}
   {} & {} & {(1,1)} & {} & {} & {} & {}  \\
   {} & {\textcolor{red}{\nearrow}} & {} & {\textcolor{red}{\searrow}} & {} & {} & {}  \\
   {(0,0)} & {} & {} & {} &{(2,0)}  & {\textcolor{green}{\longrightarrow}}  & {(3,0)}  \\
\end{array}
$$
\medskip
$$
\begin{array}{*{20}c}
   {} & {} &  {(1,1)}  & {} & {} & {} & {}  \\
   {} & {\textcolor{blue}{\nearrow}} & {} & {\textcolor{blue}{\searrow}} & {} & {} & {}  \\
   {(0,0)}  & {} & {} & {} & {(2,0)} & {\textcolor{green}{\longrightarrow}} & {(3,0)}   \\
\end{array}
$$
\normalsize
\section{The action of $\eta$}
We begin with an example

\begin{ese} \label{et}
If we consider $\mu(h,k)$ ( in this section we take ,\ without loss of generality,\  $-h$ instead of $h$) we have from (\ref{mgf})
$$u=\frac{{1 + ht - \sqrt {(1 + ht)^2  - 4kt^2 } }}{{2kt }} $$ which solved as an equation in $t$ gives
$$t=\frac{u}{ku^2-hu+1}=\sum\limits_{n=0}^{+\infty}F_n(h,k)u^{n+1} .$$
We used (\ref{genf}) with $b=h$ obtaining $F(h,k)=\left\{F_n(h,k) \right\}_{n=0}^{+\infty}$ , the generalized Fibonacci sequence.\ 
So we note that

$$\begin{cases}
\eta(F(h,k))=\mu(h,k)\\
\eta(\mu(h,k))=F(h,k)&\text{.}\\
\end{cases}$$ 

\boss
Recalling the Corollary \ref{core} we have an alternative way to find the relation proved in Corollary \ref{corL}.
When $b=h$  $$\mathcal{W}(1,h,h,k)=F(h,k)$$ and from example (\ref{et}) 
$$I^{(x)}(\mu(h,k))=\eta(\mathcal{W}(1,h-x,h-2x,x^2-hx+k)$$
$$L^{(x)}(\mu(h,k))=\eta(\mathcal{W}(1,h-x,h-x,k))=\eta(F(h-x,k))=\mu(h-x,k)$$

\eoss
\end{ese}
\noin \\The terms of $B=\eta(A)$,\ can be expressed by means of 
\noin \\ \textbf{Lagrange inversion formula \cite{Dom}}
\be
b_n=\left. {\frac{1}{{(n+1)!}}\frac{{d^{n} }}{{du^{n} }}\left\{ {\left[ {\frac{u}{{t(u)}}} \right]^{n+1} } \right\}} \right|_{u = 0}\quad. \ee

\noin Using Lagrange inversion formula we can find an analogous expression of (\ref{mom}) for $\mu_n(h,k)$. In fact
$$t(u)=\frac{u}{ku^2-hu+1}$$
thus
$$\frac{u}{{t(u)}}=ku^2-hu+1$$
and
$$\mu_n(h,k)=b_n=\left. {\frac{1}{{(n+1)!}}\frac{{d^{n} }}{{du^{n} }}\left\{ {(ku^2-hu+1)  ^{n+1} } \right\}} \right|_{u = 0}.$$
The trinomial expansion gives
$$(ku^2-hu+1)^{n+1}=
\sum\limits_{ p + q + r = n + 1 } {\frac{{(n + 1)!}}{{p!q!r!}}} ( - hu)^r (ku^2 )^q 
 $$\\
and so 
$$(ku^2-hu+1)^{n+1}=\sum\limits_{p+q=0}^{n+1} \frac{{(n + 1)!}}{{p!(n+1-p-q)!q!}}( - h)^{n+1-p-q}k^qu^{n+1-p+q}\quad .$$\\
Differentiating $n$ times
\small  
$$\frac{{d^{n} }}{{du^{n} }}\left\{ (ku^2-hu+1)^{n+1}\right\}=\sum\limits_{p+q=0}^{n+1} \frac{{(n + 1)!(n+1-p+q)!}}{{p!(n+1-p-q)!(q-p+1)!q!}}( - h)^{n+1-p-q}k^qu^{q-p+1}\quad.$$\\
\normalsize
For $u=0$ the only non zero term occurs when $q=p-1$.\ Consequently  $p+q=2p-1$ and $0\leq p+q\leq n+1$ implies $1\leq p\leq[\frac{n+2}{2}]$ then 
$$\mu_n(h,k)=\sum\limits_{p=1}^{[\frac{n+2}{2}]} \frac{{n!}}{{p!(n-2p+2)!(p-1)!}}( - h)^{n-2p+2}k^{p-1}\quad.$$\\
Taking $h=-x$ and considering odd and even values for $n$, we have
\be \label{moux} \begin{cases}
\mu_n(-x,k)=\sum\limits_{j = 0}^{\frac{{n - 1}}{2}} {A_j^{(n)} } x^{2j + 1} &\text{$n$ odd}\\
\mu_n(-x,k)=\sum\limits_{j = 0}^{\frac{{n}}{2}} {A_j^{(n)} } x^{2j} &\text{$n$ even .}\\
\end{cases}\ee
where
\be \label{amoux}
A_j^{(n)}=\begin{cases}
\frac{1}{{n + 1}}\left( {\begin{array}{*{20}c}
   {n + 1}  \\
   {\begin{array}{*{20}c}
   {\frac{{n + 1}}{2} - j,} & {2j + 1,} & {\frac{{n - 1}}{2} - j}  \\
\end{array}}  \\
\end{array}} \right)k^{\frac{{n - 1}}{2} - j} &\text{$n$ odd}\\
\frac{1}{{n + 1}}\left( {\begin{array}{*{20}c}
   {n + 1}  \\
   {\begin{array}{*{20}c}
   {\frac{n}{2} + 1 - j,} & {2j,} & {\frac{n}{2} - j}  \\
\end{array}}  \\
\end{array}} \right)k^{\frac{n}{2} - j} &\text{$n$ even .}\\
\end{cases}
\ee
\boss
\itshape 
\noin\\
If $x=0$ for $n=2m$ we have   $$\mu_{2m}(0,k)=k^mC_m \quad\quad C_m= 
\frac{1}{{2m + 1}}\binom{2m+1}{m}\quad $$ 
and $C_m$ is the m-th Catalan number, while if $n=2m+1$ we have $$\mu_{2m+1}(0,k)=0\quad.$$ 
\eoss
\boss[\textbf{Orthogonality relations}]
We consider the polynomial $P_n(x)=\mathcal{W}(1,h+x,h+x,k)$.\ Its explicit expression can be found observing that from (\ref{genf})
we have
\scriptsize
$$\sum\limits_{n = 0}^{ + \infty } {P_n (x)t^n }  = \frac{1}{{1 - (h + x)t + kt^2 }}= \sum\limits_{j = 0}^{ + \infty } {((h + x)t - kt^2 )^j }  = \sum\limits_{j = 0}^{ + \infty } \sum\limits_{l = 0}^j \binom{j}{l}(h + x)^{j - l} ( - k)^l t^{l + j}\quad. $$
\normalsize
Rearranging indexes and posing  $l+j=n$ we obtain
\be \label{pienne}
P_n (x) = \sum\limits_{l = 0}^{\left\lfloor  {\frac{n}{2}}\right\rfloor}\binom{n - l}{l}(h + x)^{n - 2l} ( - k)^l .\ee
and  the generic coefficient $P^{(n)}_j$ of $x^j$ follows from the $j$-th derivative: 
\be \label{piennej}
P_j^{(n)}  = \sum\limits_{l = 0}^{\left\lfloor{\frac{n}{2}}\right\rfloor}\binom{n - l}{l} \binom{n - 2l}{j}h^{n - 2l - j} ( - k)^l \quad.\ee
Now from definition of the functional $\mathcal{V}$ \cite{Chi} 
\be \label{orto}
\left\{ \begin{array}{l}
\mathcal{V}[1]=1\\ 
\mathcal{V}[P_m(x)P_n(x)]=0 \ \text{if $m\neq n$}\\
\end{array} \right.
 \ee
from (\ref{pol}) $P_0(x)=1$ and from (\ref{orto}) choosing $m=0$ we have
$$\mathcal{V}[P_n(x)]= \delta(n,0)$$
which in this case becomes the following relation 
\scriptsize
$$
\sum\limits_{j = 0}^n P_j^{(n)} \mu _j (h,k) = \sum\limits_{j = 0}^n {\sum\limits_{l = 0}^{ \left\lfloor{\frac{n}{2}}\right\rfloor} \binom{n - l}{l} \binom{n - 2l}{j}h^{n - 2l - j} ( - k)^l }\sum\limits_{p = 1}^{\left\lfloor {\frac{{j + 2}}{2}} \right\rfloor} {\frac{{j!( - h)^{j - 2p + 2} k^{p - 1}}}{{p!(j - 2p + 2)!(p - 1)!}}}= \delta (n,0)\quad .$$
\normalsize
And when $h=0$ we have
\be\label{orth1} P_n(x)=\sum\limits_{i = 0}^{\left\lfloor {\frac{n}{2}} \right\rfloor}\binom{n - i}{i}( - k)^i x^{n - 2i}=E_n(x,k)\quad . \ee
\eoss
\noin where $E_n(x,k)$ is the $n$-th Dickson polynomial of the second kind \cite{Lidl} .

\noin \\ So if $n=2m$, recalling that $\mu_{2m}(0,k)=k^m C_m$, we obtain a similar orthogonality relation where Catalan numbers are involved
\be \label{catid}
\sum\limits_{i = 0}^{m}\binom{2m - i}{i}( - k)^i k^{m - i} C_{m - i}= \delta (m,0)\quad. 
\ee

\medskip

\noin As far as we know this Catalan identity is new. Of course (\ref{catid}) is not difficult to prove (try it with Zeilberger's program \cite{PWZ}, 
for example), but it seems interesting also for the context it rises from.

\end{document}